%% file: agt-5-58.tex
\documentclass{gtart_h}
\input agtout

\lognumber{58}
\volumenumber{5}
\volumeyear{2005}
\papernumber{58}
\pagenumbers{1471}{1479}
\received{9 August 2005} 
\accepted{20 September 2005}
\published{30 October 2005}

\usepackage{amsmath,amssymb}
\usepackage{graphicx,psfrag}

\def\psfraga <#1,#2> #3#4{%
\psfrag {#3}{\smash{\rlap{\kern #1 \raise #2\hbox{#4}}}}}

\def\figref#1{\hyperlink{#1anchor}{Figure~\ref*{#1}}}
\def\anchor#1{\noindent\hypertarget{#1anchor}{\smash{$\phantom{99}$}}\newline}

\newcommand{\R}{\mathbb R}
\newcommand{\Z}{\mathbb Z}
\newcommand{\onen}{1,\dots,n}
\newcommand{\ioner}{{i_1\dots i_r}}
\newcommand{\Gd}{\delta}

\newcommand{\Ge}{\varepsilon}

\newcommand{\Gt}{\theta}
\newcommand{\mb}{\bar{\mu}}
\newcommand{\tp}{\tilde\pi}
\newcommand{\tF}{\widetilde{F}}
\newcommand{\tZ}{\widetilde{Z}}

\newtheorem{thm}{Theorem}[section]
\newtheorem{lem}[thm]{Lemma}

\theoremstyle{definition}
\newtheorem{ex}[thm]{Example}
\theoremstyle{definition}

\theoremstyle{remark}

\newcommand{\lk}{\operatorname{lk}}
\def\<{\langle}
\def\>{\rangle}

\begin{document}

\title{Skein relations for Milnor's $\mu$-invariants}
\author{Michael Polyak}
\address{Department of mathematics, Technion, Haifa 32000, Israel}
\email{polyak@math.technion.ac.il}

\begin{abstract}
The theory of link-homotopy, introduced by Milnor, is an important
part of the knot theory, with Milnor's $\mb$-invariants being the
basic set of link-homotopy invariants.
Skein relations for knot and link invariants played a crucial role
in the recent developments of knot theory.
However, while skein relations for Alexander and Jones invariants
are known for quite a while, a similar treatment of Milnor's
$\mb$-invariants was missing.
We fill this gap by deducing simple skein relations for
link-homotopy $\mu$-invariants of string links.
\end{abstract}
\asciiabstract{%
The theory of link-homotopy, introduced by Milnor, is an important
part of the knot theory, with Milnor's mu-bar-invariants being the
basic set of link-homotopy invariants.  Skein relations for knot and
link invariants played a crucial role in the recent developments of
knot theory.  However, while skein relations for Alexander and Jones
invariants are known for quite a while, a similar treatment of
Milnor's mu-bar-invariants was missing.  We fill this gap by deducing
simple skein relations for link-homotopy mu-invariants of string
links.}



\asciikeywords{String links, link homotopy, Milnor's mu-invariants,
skein relations}

\maketitle

\section{Introduction}

\subsection{Some history}
Skein relations for various invariants played a significant role in
knot theory and related areas in the last decade.
Probably the first widely-used skein relation, discovered by Conway,
was the skein relation for the Conway-Alexander polynomial.
Apart from multiple technical consequences, it led to the discovery
of similar skein relations for other knot invariants, notably for
the Jones and HOMFLY polynomials.
In the theory of 3-manifolds, the skein relation for the Casson
invariant of homology spheres was also fruitfully used in different
contexts.

An important part of the knot theory is the theory of link-homotopy,
initiated by Milnor.
Link-homotopy is a useful notion to isolate the linking fenomena
from the self-knotting ones and to study it separately.
A celebrated example of link-homotopy invariants is given by
Milnor's $\mb_{\ioner,j}$ invariants \cite{Mi,Mi2} with non-repeating
indices $1\le i_1,\dots i_r,j\le n$.
Roughly speaking, these describe the dependence of $j$-th parallel
on the meridians of $i_1$-th, \dots, $i_r$-th components.
The simplest invariant $\mb_{i,j}$ is just the linking number of the
corresponding components.
The next one, $\mb_{i_1i_2,j}$, detects the Borromean-type linking
of the corresponding 3 components and, together with the linking
numbers, classify 3-component links up to link-homotopy.

Unfortunately, a complicated self-recurrent indeterminacy in the
definition of $\mb$-invariants (reflected in the use of notation
$\mb$, rather than $\mu$) for a long time slowed down their study.
The introduction of string links \cite{HL} considerably improved
the situation, since a version of $\mb$-invariants modified for
string links is free of this indeterminacy; thus (and to stress a
special role of the $j$-th component) we will further use the
notation $\mu_{i_1\dots i_r}(L, L_j)$ for these invariants.
Milnor's invariants classify string links up to link-homotopy
\cite{HL}. Surprisingly, up to now $\mu$-invariants remained aside
from the well-developed scheme of skein relations.

\subsection{Brief statement of results}
We deduce new skein relations for $\mu$-invariants. Usually in the
knot theory skein relations involve links, obtained by different
splittings of a diagram in a crossing. In the context of string
links this leads to an appearance of tangles which are not pure
any more, but contain a new ``loose" component. Fortunately, it is
easy to extend the definition of $\mu$-invariants to such tangles.
We next note that $\mu_{i_1\dots i_r}(L, l)=0$ for any string link
$L$ whose loose component $l$ passes everywhere in front of the
other strings. Thus it suffices to study the jump of
$\mu_{i_1\dots i_r}$ under a crossing change of $l$ with any
other, say, $i_k$-th, component. It turns out, that this jump can
be expressed via the invariants $\mu_{i_1\dots
i_{k-1}}(L-L_{i_k},l_0)$ and $\mu_{i_{k+1}\dots
i_r}(L-L_{i_k},l_\infty)$ of string links with new loose
components $l_0$ and $l_\infty$ obtained by splitting the crossing
in two possible ways, see Section \ref{main_sec}.

This research was partially supported by C. Wellner Research fund
and 6th EC Programme ``Structuring the European Research Area",
contract Nr.\ RITA-CT-2004-505493

\section{Preliminaries}

\subsection{String links and link-homotopy}
\begin{figure}[htb]\anchor{string_fig}\small
\psfraga <-2pt,0pt> {0}{$0$}
\psfraga <-2pt,0pt> {1}{$1$}
\psfrag {x}{$x$}
\psfraga <-2pt,0pt> {y}{$y$}
\psfraga <-2pt,0pt> {z}{$z$}
\psfraga <-2pt,0pt> {m1}{$m_1$}
\psfraga <-2pt,0pt> {m2}{$m_2$}
\psfraga <-2pt,0pt> {m3}{$m_3$}
\centerline{\includegraphics[width=4.7in]{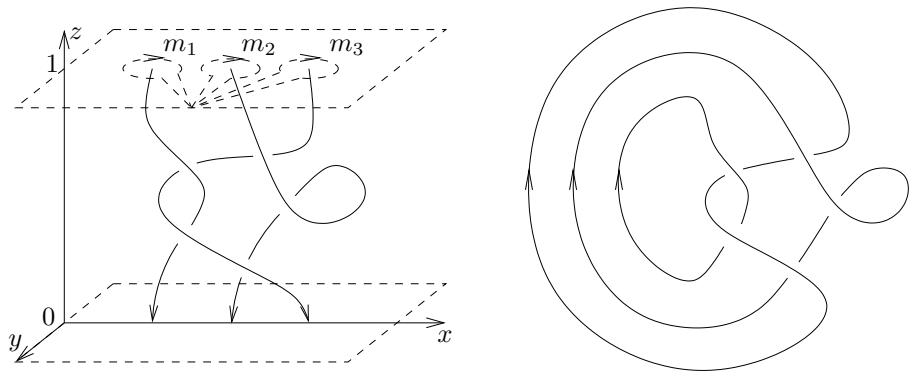}}
\caption{\label{string_fig} String link and its closure}
\end{figure}

Let $D^2$ be a disc in the $xy$-plane which intersects the
$x$-axis. An {\em $n$-component string link} (see \cite{Le},
\cite{HL}) $L$ is an ordered collection of $n$ disjoint arcs
properly embedded in $D^2\times[0,1]$ in such a way, that the
$i$-th arc ends in the points $p_i\times\{0,1\}$, where
$p_i=(x_i,0)$ are some prescribed points on the $x$-axis,
enumerated in the natural order $x_1<x_2<\dots<x_n$. We assume
that all arcs of $L$ are oriented downwards. By the {\em closure}
$\overline{L}$ of a string link $L$ we mean the braid closure of
$L$. It is an $n$-component link obtained from $L$ by an addition
of $n$ disjoint arcs in the plane $\{y=0\}$, each of which meets
$D^2\times[0,1]$ only at the endpoints $p_i\times\{0,1\}$ of $L$,
as illustrated in \figref{string_fig}. The linking number
$\lk$ of two components of $L$ is their linking number in
$\overline{L}$.
Two string links are {\em link-homotopic}, if one can be transformed
into the other by homotopy, which fails to be isotopy only in a
finite number of instants, when a (generic) self-intersection point
appears on one of the arcs.

\subsection{String links with a loose component}
\begin{figure}[htb]\anchor{loose_fig}\small
\psfraga <-2pt,0pt> {L}{$L$}
\psfraga <-2pt,0pt> {l}{$l$}
\psfraga <-2pt,0pt> {Ll}{$(L,l)$}
\centerline{\includegraphics[width=4.9in]{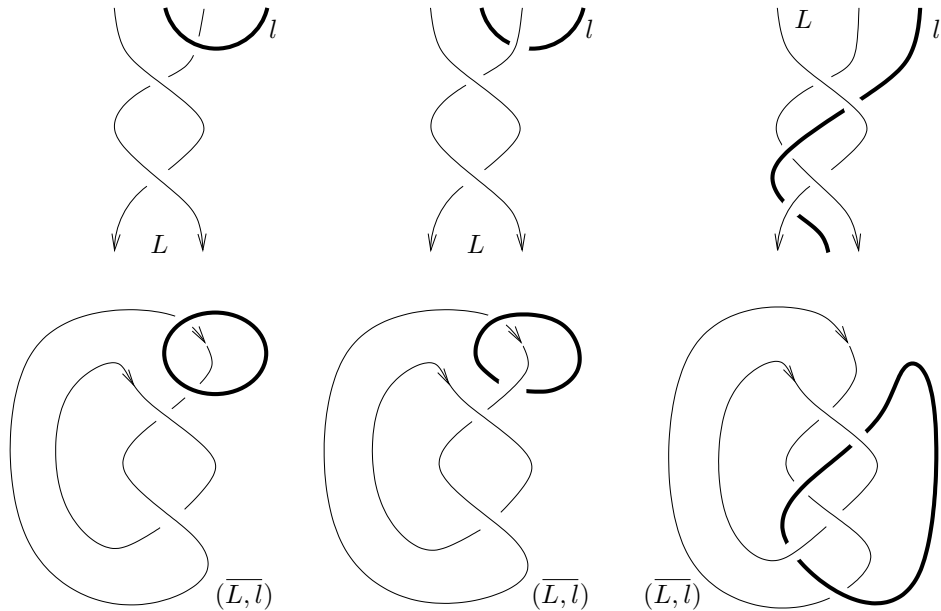}}
\caption{\label{loose_fig} 2-component string links with a loose
component and their closures}
\end{figure}

Further we will consider a more general class of tangles. A {\em
$n$-component string link $L$ with an additional loose component
$l$} is a pair $(L, l)$ which consists of an $n$-string link $L$
together with an additional oriented arc $l$ properly embedded in
$D^2\times[0,1]$ in such a way, that $l\cap L=\emptyset$ and $l$
starts on $D^2\times\{1\}$ and ends on $D^2\times\{0,1\}$. See
\figref{loose_fig}, where loose components are depicted in
bold. By the {\em closure} $\overline{(L,l)}$ of a string link $L$
with a loose component $l$ we mean the closure of $L$, together
with a closure of $l$ by a standard arc in $\R^3-D^2\times(0,1)$,
which meets $l$ only at its endpoints and passes ``in front" of
$L$ i.e.\ lies in the half-space $y\ge0$, as illustrated in 
\figref{loose_fig}.

In particular, for any $n$-string link $L=\cup_{i=1}^n L_i$ and
$1\ge j\ge n$ one may consider $(L-L_j,L_j)$ as an $(n-1)$-string
link with an additional loose component $L_j$. In this case it
does not matter whether to close $L_j$ in front of $L$ or in the
plane $\{y=0\}$ and closures $\overline{L}$ and
$\overline{(L-L_j,L_j)}$ give isotopic links.

\subsection{Magnus expansion}
Milnor's link-homotopy $\mu$-invariants \cite{Mi} can be defined
in several ways. We choose a construction most suitable for our
purposes and refer the reader to Milnor's work \cite{Mi,Mi2} and its
adaptations to string links (e.g.\ \cite{Le}, \cite{HL}) for the
general case.

Let $L=\cup_{i=1}^nL_i$ be an $n$-component string link. Consider
the link group $\pi=\pi_1(D^2\times[0,1]-L)$ (with the base point
on the upper boundary disc $D^2\times\{1\}$). Denote by
$m_i\in\pi$, $i=\onen$ the {\em canonical meridians} represented
by the standard non-intersecting curves in $D^2\times\{1\}$ with
$\lk(m_i,L_i)=+1$, as shown in \figref{string_fig}. If $L$
would be a braid, these meridians would freely generate $\pi$,
with any other meridian of $L_i$ in $\pi$ being a conjugate of
$m_i$. For the string links, similar results hold for the reduced
link group $\tp$.

For any group $G$ with a finite set of generators $x_1,\dots x_n$,
the {\em reduced group} $\tilde{G}$ is the factor group of $G$ by
relations $[x'_i,x''_i]=1$, $i=\onen$ for any two elements $x'_i$
and $x''_i$ in the conjugacy class of $x_i$.

Proceeding similarly to the usual construction of Wirtinger's
presentation, one can show (see \cite{HL}) that $\tp$ is generated
by $m_i$. Let $F$ be the free group on $n$ generators $x_1,\dots
x_n$. The map $F\to\pi$ defined by $x_i\mapsto m_i$ induces the
isomorphism $\tF\cong\tp$ of the reduced groups \cite{HL}. We will
use the same notation for the elements of $\pi$ and their images
in $\tp\cong\tF$.

Now, let $\Z[[X_1,\dots,X_n]]$ be the ring of power series in $n$
non-commuting variables $X_i$ and denote by $\tZ$ its factor by
all the monomials, where at least one of the generators appears
more than once. The {\em Magnus expansion} is a ring homomorphism
of the group ring $\Z F$ into $\Z[[X_1,\dots,X_n]]$, defined by
$x_i\mapsto 1+X_i$. It induces the homomorphisms $\Gt:\Z\tF\to\tZ$
and $\Gt_L:\Z\tp\to\tZ$ of the corresponding reduced group rings.

\subsection{Milnor's $\mu$-invariants}
Let $l$ an immersed arc in $D^2\times[0,1]-L$ with $\partial l\in
R^1\times\{0\}\times\{0,1\}$. Its closure by an arc on the
boundary of the cylinder $D^2\times[0,1]$ in front of $L$ (i.e.\ in
the half-space $y\ge0$) gives a well-defined loop $\bar{l}$ in
$\tp$. Thus one may define the {\em Milnor's invariants}
$\mu_{\ioner}(L,l)$ of a string link with a loose component
$(L,l)$ as coefficients of the Magnus expansion $\Gt_L(\bar{l})$
of $\bar{l}$:
$$\Gt_L(l)=\sum\mu_{\ioner}(L,l)X_{i_1}X_{i_2}\dots X_{i_r}\ .$$
Note that if $l$ goes in front of $L$, i.e.\ overpasses all other
components on the $xz$-plane diagram of $L\cup l$, then all
invariants $\mu_{\ioner}(L,l)$ vanish. Indeed, we choose to close
$l$ in the front half-space
$\{y\ge0\}\cap\partial(D^2\times[0,1])$, so the loop $\bar{l}$ is
trivial in $\pi$ (and hence trivial in $\tp$).

For string links, Milnor's invariants $\mu_{\ioner}(L-L_j,L_j)$
are usually denoted by $\mu_{\ioner,j}(L)$ (see e.g.\ \cite{Le},
\cite{HL}). Modulo lower degree invariants
$\mu_{\ioner}(L-L_j,L_j)\equiv\mb_{\ioner,j}(\overline{L})$, where
$\mb_{\ioner,j}(\overline{L})$ are the original Milnor's link
invariants \cite{Mi,Mi2}.

\subsection{Fox's free calculus}
Instead of the Magnus expansion, $\mu$-invariants may be defined
via Fox free calculus.

Fox's {\em free derivatives} $\Gd_i:\Z F\to \Z F$, $i=\onen$, are
defined by putting $\Gd_i1=0$, $\Gd_ix_i=1$, $\Gd_ix_j=0$, $j\ne i$
and extending the function $\Gd_i$ to $\Z F$ by linearity and the
rule $\Gd_i(uv)=\Gd_iu+u\cdot\Gd_iv$, $u,v\in Z F$.
In particular, it is easy to see that $\Gd_i(x_i^{-1})=-x_i$.
For $u\in\Z F$, denote by $\Gd_{i_1}\dots\Gd_{i_r}u(1)\in\Z$ the
value of $\Gd_{i_1}\dots\Gd_{i_r}u$ in $x_1=\dots=x_n=1$.

Free derivatives in $\Z F$ induce similar Fox derivatives
$\Gd_i:\Z\tF\to\Z\tF$ and $\Gd_i:\Z\tp\to\Z\tp$ in the reduced
group rings. From the definition of $\mu$-invariants we conclude
that $\mu_{\ioner}(L,l)=\Gd_{i_1}\dots\Gd_{i_r}l(1)$.

\section{Skein relations for $\mu$-invariants}

\subsection{Main Theorem}\label{main_sec}
\begin{figure}[htb]\anchor{split_fig}\small
\psfraga <-2pt,0pt> {L}{$L_{i_k}$}
\psfraga <-2pt,0pt> {l-}{$l_-$}
\psfraga <-2pt,0pt> {l+}{$l_+$}
\psfraga <-2pt,0pt> {l0}{$l_0$}
\psfraga <0pt,0pt> {l}{$l$}
\centerline{\includegraphics[width=4.7in]{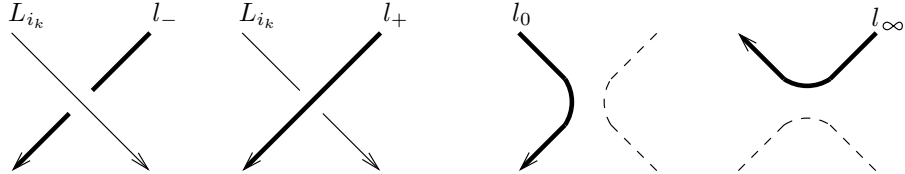}}
\caption{\label{split_fig} Splitting the link in a crossing}
\end{figure}

Consider two $n$-string links $(L,l_+)$ and $(L,l_-)$ with loose
components such that their diagrams coincide everywhere, except
for a crossing $d$, where $l_+$ has the positive crossing and
$l_-$ the negative crossing with $i_k$-th component $L_{i_k}$ of
$L$, see \figref{split_fig}.
We define two new $(n-1)$-string links with loose components by
splitting the diagram in $d$ in two possible ways. Let $l_\infty$
be the component of the splitting going along $l$, respecting the
orientation, and then switching to $L_{i_k}$ against the
orientation. Let also $l_0$ be the component going first along
$L_{i_k}$ and then along $l$, respecting the orientation. Each of
$(L-L_{i_k},l_0)$ and $(L-L_{i_k},l_\infty)$ is a $(n-1)$-string
link with a loose component, see \figref{split_fig}.

\begin{thm}\label{main_th}
Let $L$, $l_\pm$, $l_0$ and $l_\infty$ be as above. Then
$$\mu_{i_1\dots i_k\dots i_r}(L,l_+)-
\mu_{i_1\dots i_k\dots i_r}(L,l_-)= \mu_{i_1\dots
i_{k-1}}(L-L_{i_k},l_\infty) \cdot\mu_{i_{k+1}\dots
i_r}(L-L_{i_k},l_0).$$
Here it is understood that in particular
cases $k=1$ or $k=r$ we have
$$\mu_{i_k}(L,l_+)-\mu_{i_k}(L,l_-)=1;$$
$$\mu_{i_1\dots i_k}(L,l_+)-\mu_{i_1\dots i_k}(L,l_-)=
\mu_{i_1\dots i_{k-1}}(L-L_{i_k},l_\infty);$$
$$\mu_{i_k\dots i_r}(L,l_+)-\mu_{i_k\dots i_r}(L,l_-)=
\mu_{i_{k+1}\dots i_r}(L-L_{i_k},l_0).$$
\end{thm}

\begin{figure}[htb]\anchor{example_fig}\small
\psfraga <-2pt,0pt> {l=}{$l=l_-$}
\psfraga <-2pt,0pt> {l+}{$l_+$}
\psfraga <-2pt,0pt> {l0}{$l_0$}
\psfraga <0pt,0pt> {l}{$l$}
\centerline{
\includegraphics[height=1.3in,width=4.6in]{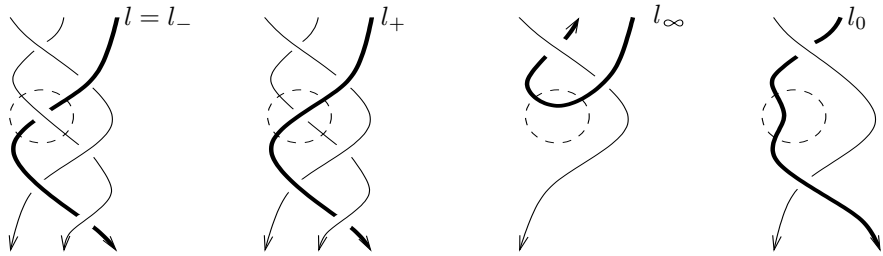}}
\caption{\label{example_fig} A computation of $\mu_{12,3}$ for
Borromean rings}
\end{figure}

\begin{ex}
Consider the string link $(L,l)$ depicted in 
\figref{example_fig} and let us compute $\mu_{12}(L,l)$. Notice that
if we switch the crossing $d$ to the positive one, we get a link
$(L,l_+)$ with $L_1$ unlinked from $L_2$ and $l_+$, so
$\mu_{12}(L,l_+)=0$. Thus
$\mu_{12}(L,l)=\mu_{12}(L,l)-\mu_{12}(L,l_+)=
-\mu_{1}(L_1,l_\infty)=-1$, in agreement with the fact that the
closure of $(L,l)$ is the (negative) Borromean link.
\end{ex}

\subsection{The proof of Theorem \ref{main_th}}
We will need the following simple key fact about Fox free
calculus, which we leave to the reader as an exercise:

\begin{lem}\label{fox_lem}
Let $1\le k\le r<j\le n$ and $\Ge=\pm1$.
Then for any $u,v\in F$
\begin{equation}\label{fox1_eq}
\Gd_1\dots\Gd_r(ux_k^\Ge v)(1)=
\Ge\cdot\Gd_1\dots\Gd_{k-1}u(1)\cdot\Gd_{k+1}\dots\Gd_rv(1)+
\Gd_1\dots\Gd_r(uv)(1).
\end{equation}
\end{lem}


Also, to simplify the notations in the statement of Theorem
\ref{main_th}, let us reorder the components so that $i_m\to m$,
$m=1,\dots r$ (and $j\ne1,\dots,r$).

\begin{proof}[Proof of Theorem \ref{main_th}]
The component $l_+$ passes either over, or under $L_k$ in the
crossing $d$. Let us consider these cases separately. Assume first
that $l_+$ passes under $L_k$ in $d$. Then $l_+=vu^{-1}x_kuw$,
where $v$ and $w$ are the parts of $l_+$ before and after $d$, and
$u$ is the part of $l_k$ before $d$, see \figref{proof1_fig}.
\begin{figure}[htb]\anchor{proof1_fig}\small
\psfraga <-2pt,0pt> {u}{$u$}
\psfraga <-2pt,0pt> {v}{$v$}
\psfraga <-2pt,0pt> {w}{$w$}
\psfraga <-2pt,0pt> {Lk}{$L_k$}
\psfraga <-2pt,0pt> {l+}{$l_+$}
\psfraga <-2pt,0pt> {l-}{$l_-$}
\psfraga <-2pt,0pt> {ux}{$u^{-1}x_ku$}
\centerline{\includegraphics[width=4.7in]{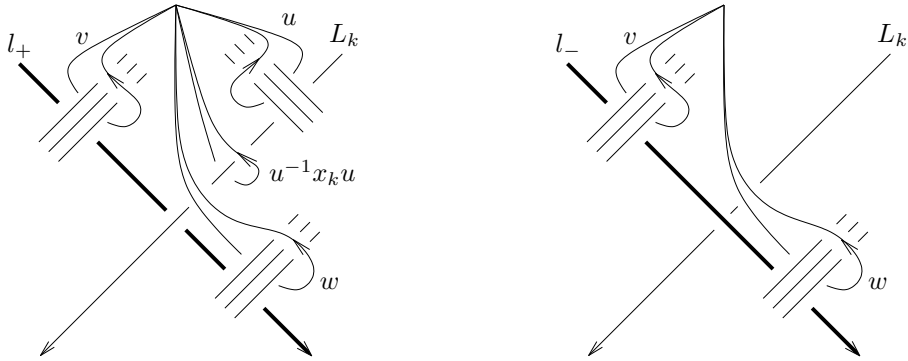}}
\caption{\label{proof1_fig} Computing $l_+$ and $l_-$}
\end{figure}
Thus, by the definition of $\mu_{1\dots k\dots r}(L,l)$ and
\eqref{fox1_eq} we have
\begin{multline*}
\mu_{1\dots r}(L,l_+)=\Gd_1\dots\Gd_r(vu^{-1}x_kuw)(1)=\\
\Gd_1\dots\Gd_{k-1}(vu^{-1})(1)\cdot\Gd_{k+1}\dots\Gd_r(uw)(1)+
\Gd_1\dots\Gd_r(vw)(1).
\end{multline*}
Now, $l_-$ passes over $L_k$, so $l_-=vw$, where again $v$ and $w$
are the parts of $l_-$ before and after $d$, see 
\figref{proof1_fig}. Thus
$$\mu_{1\dots k\dots r}(L,l_-)=\Gd_1\dots\Gd_k\dots\Gd_r(vw)(1).$$
Therefore
$$\mu_{1\dots k\dots r}(L,l_+)-\mu_{1\dots k\dots r}(L,l_-)=
\Gd_1\dots\Gd_{k-1}(vu^{-1})(1)\cdot\Gd_{k+1}\dots\Gd_r(uw)(1).$$
But $vu^{-1}$ and $uw$ are exactly $l_\infty$ and $l_0$, see
\figref{proof2_fig}.

\begin{figure}[ht!]\anchor{proof2_fig}\small
\psfraga <-2pt,0pt> {u}{$u$}
\psfraga <-2pt,0pt> {u-}{$u^{-1}$}
\psfraga <-2pt,0pt> {v}{$v$}
\psfraga <-2pt,0pt> {w}{$w$}
\psfraga <-2pt,0pt> {l0}{$l_0$}
\psfraga <-0pt,0pt> {l}{$l$}
\centerline{\includegraphics[width=4.7in]{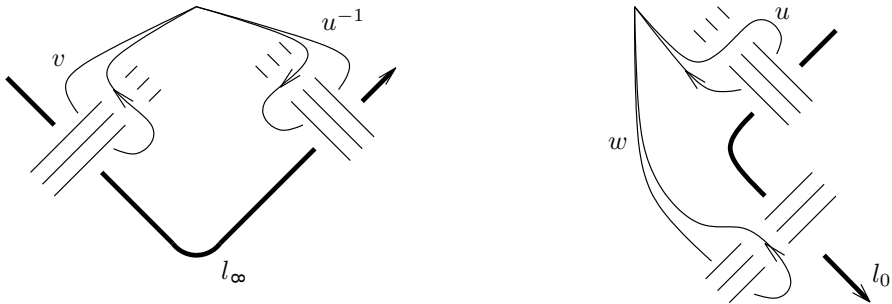}}
\caption{\label{proof2_fig} Computing $l_\infty$ and $l_0$}
\end{figure}

It remains to use again the definition of $\mu$-invariants to
obtain
$$\mu_{1\dots k\dots r}(L,l_+)-\mu_{1\dots k\dots r}(L,l_-)=
\mu_{1\dots k-1}(L,l_\infty)\cdot\mu_{k+1\dots r}(L,l_0).$$

In the case when $l_+$ passes over $L_k$ at the crossing $d$, the
proof is completely similar. This time $l_-=vu^{-1}x_k^{-1}uw$ and
$l_+=vw$, thus basically we just exchange $l_+$ and $l_-$ in the
previous proof. Also, the degree of $x_k$ switches from $+1$ to
$-1$, which leads to the negative sign in the application of
\eqref{fox1_eq} and cancels out with the sign appearing from the
exchange of $L_+$ and $L_-$.
\end{proof}

\Addresses\recd
\end{document}

%% file: agtout.tex
\def\ifplaintex{\expandafter\ifx\csname documentclass\endcsname\relax}

\def\gtp{{\mathsurround=0pt\it $\cal G\mskip-2mu$eometry \&\ 
$\cal T\!\!$opology $\cal P\!$ublications}}  

\def\recd{{\small Received:\qua\receiveddate\ifx\reviseddate\relax
\else\qquad Revised:\qua\reviseddate\fi\par}} 

\def\lognumber#1{\def\thelognumber{#1}}
\def\volumenumber#1{\def\thevolumenumber{#1}}
\def\volumeyear#1{\def\thevolumeyear{#1}}
\def\papernumber#1{\def\thepapernumber{#1}}
\def\pagenumbers#1#2{\def\startpage{#1}\def\finishpage{#2}}
\def\published#1{\def\publishdate{#1}}

\def\received#1{\def\receiveddate{#1}}

\def\accepted#1{\def\accepteddate{#1}}

\long\def\asciiabstract#1{\long\def\theasciiabstract{#1}}
\def\asciikeywords#1{\def\theasciikeywords{#1}}

\let\\\par\let\thelognumber\relax\let\thevolumenumber\relax
\let\thepapernumber\relax\let\thevolumeyear\relax\let\startpage\relax
\let\finishpage\relax\let\publishdate\relax\let\receiveddate\relax
\let\reviseddate\relax\let\accepteddate\relax\let\theasciititle\relax
\let\theasciiauthors\relax
\let\theasciiabstract\relax\let\theasciikeywords\relax

\let\theasciiemail\relax

\ifplaintex
\font\logobig=cmssbx10 scaled 3836
\font\logomed=cmssbx10 scaled 2557
\else
\font\logobig=cmssbx10 scaled 4200
\font\logomed=cmssbx10 scaled 2800
\fi

\long\def\makeagttitle{   
\count0=\startpage
\agt\hfill      
\hbox to 45truept{\vbox to 0pt{\vglue -13truept{\logomed A\kern -.37em{\logobig 
T}\kern -.38em G}\vss}\hss}
\break
{\small Volume \thevolumenumber\ (\thevolumeyear)
\startpage--\finishpage\nl
Published: \publishdate}

\vglue .25truein

{\parskip=0pt\leftskip 0pt plus
1fil\def\\{\par\smallskip}{\Large\bf\thetitle}\par\medskip} \vglue
0.05truein

{\parskip=0pt\leftskip 0pt plus 1fil\def\\{\par}{\sc\theauthors}
\par\medskip}%
 
\vglue 0.03truein 

{\small\leftskip 25truept\rightskip 25truept{\bf Abstract}\stdspace\theabstract

{\bf AMS Classification}\stdspace\theprimaryclass
\ifx\thesecondaryclass\relax\else; \thesecondaryclass\fi\par
{\bf Keywords}\stdspace \thekeywords\par}\vglue 7truept

}   

\ifplaintex
\font\phead=cmsl9 scaled 950
\font\pnum=cmbx10 scaled 913
\font\pfoot=cmsl9 scaled 950
\headline{\vbox to 0pt{\vskip -4.5mm\line{\small\phead\ifnum
\count0=\startpage ISSN 1472-2739 (on-line) 1472-2747 (printed)
\hfill {\pnum\folio}\else\ifodd\count0\def\\{ }%
\ifx\theshorttitle\relax\thetitle\else\theshorttitle\fi\hfill{\pnum\folio}
\else\def\\{ and }{\pnum\folio}\hfill\ifx\theshortauthors\relax\theauthors
\else\theshortauthors\fi\fi\fi}\vss}}
\footline{\vbox to 0pt{\vglue 0mm\line{\small\pfoot\ifnum\count0=\startpage
\copyright\ \gtp\hfill\else
\agt, Volume \thevolumenumber\ (\thevolumeyear)\hfill\fi}\vss}}
\else
\font=cmsl9 scaled 1050
\font\lnum=cmbx10 
\font=cmsl9 scaled 1050
\makeatletter
\def\@oddhead{{\small\lhead\ifnum\count0=\startpage ISSN 1472-2739 
(on-line) 1472-2747 (printed)\hfill {\lnum\number\count0}\else\ifodd\count0
\def\\{ }\ifx\theshorttitle\relax \thetitle \else\theshorttitle\fi\hfill
{\lnum\number\count0}\else\def\\{ and }{\lnum\number\count0}
\hfill\ifx\theshortauthors\relax 
\theauthors\else\theshortauthors\fi\fi\fi}}\def\@evenhead{\@oddhead}
\def\@oddfoot{\small\lfoot\ifnum\count0=\startpage\copyright\ \gtp\hfill\else
\agt, Volume \thevolumenumber\ (\thevolumeyear)\hfill\fi}
\def\@evenfoot{\@oddfoot}
\makeatother
\fi
\let\maketitlepage\makeagttitle
\let\makeshorttitle\maketitlepage
\let\maketitle\maketitlepage

\newwrite\gtoutfile
\long\gdef\makeheadfile{  
{\def\\{, }\def\s{ }
\immediate\openout\gtoutfile head.xxx
\immediate\write\gtoutfile{Proxy-for: \ifx\theasciiauthors\relax
\theauthors\else\theasciiauthors\fi\s<\ifx\theasciiemail\relax\theemail\else\theasciiemail\fi>}
\immediate\write\gtoutfile{\noexpand\\}
\immediate\write\gtoutfile{Authors: \ifx\theasciiauthors\relax
\theauthors\else\theasciiauthors\fi}
{\def\\{ }\immediate\write\gtoutfile{Title: \ifx\theasciititle\relax
\thetitle\else\theasciititle\fi}}
\immediate\write\gtoutfile{Subj-class: GT or SG, GR etc}
\immediate\write\gtoutfile{MSC-class: \theprimaryclass\ifx\thesecondaryclass\relax\else, \thesecondaryclass\fi}
\immediate\write\gtoutfile{Journal-ref: Algebr. Geom. Topol. \thevolumenumber\s
(\thevolumeyear) \startpage-\finishpage}
\immediate\write\gtoutfile{Comments: Published by Algebraic and
Geometric Topology at}
\immediate\write\gtoutfile{\s\s\s  http://www.maths.warwick.ac.uk/agt/AGTVol\thevolumenumber/agt-\thevolumenumber-\thepapernumber.abs.html}
\immediate\write\gtoutfile{\noexpand\\}
\immediate\write\gtoutfile{}
\ifx\theasciiabstract\relax
\immediate\write\gtoutfile{\theabstract}\else
\immediate\write\gtoutfile{\theasciiabstract}\fi
\immediate\write\gtoutfile{}
\immediate\write\gtoutfile{\noexpand\\}
\immediate\write\gtoutfile{}
\immediate\closeout\gtoutfile}}  

\def\maketitlepage{\makeagttitle\makeheadfile}
\let\makeshorttitle\maketitlepage
\let\maketitle\maketitlepage